\documentclass[11pt]{article}
\usepackage{amssymb,epsfig,amsmath}
\usepackage[dvips]{color}
\usepackage[T1]{fontenc}
\usepackage[english]{babel}
\usepackage[utf8]{inputenc}
\usepackage{comment}
\usepackage{hyperref}
\usepackage{graphicx}
\usepackage{multicol}
\usepackage{caption}
\newtheorem{defi}{Definition}
\newtheorem{prop}[defi]{Proposition}
\newtheorem{theo}[defi]{Theorem}
\newtheorem{conj}[defi]{Conjecture}
\newtheorem{lemm}[defi]{Lemma}
\newtheorem{coro}[defi]{Corollary}
\newtheorem{rema}[defi]{Remark}
\newtheorem{exem}[defi]{Example}
\newtheorem{exems}[defi]{Examples}

\newcommand{\bdefi}{\begin{defi}}
\newcommand{\edefi}{\end{defi}}
\newcommand{\bprop}{\begin{prop}}
\newcommand{\eprop}{\end{prop}}
\newcommand{\btheo}{\begin{theo}}
\newcommand{\etheo}{\end{theo}}
\newcommand{\blemm}{\begin{lemm}}
\newcommand{\brema}{\begin{rema}}
\newcommand{\erema}{\end{rema}}
\newcommand{\bexer}{\begin{exem}}
\newcommand{\eexer}{\end{exem}}
\newcommand{\bexems}{\begin{exems}}
\newcommand{\eexems}{\end{exems}}
\newcommand{\bconj}{\begin{conj}}
\newcommand{\econj}{\end{conj}}
\newcommand{\elemm}{\end{lemm}}
\newcommand{\bcoro}{\begin{coro}}
\newcommand{\ecoro}{\end{coro}}

\newcommand{\rem}{\noindent{\bf Remark. }}

\usepackage{mathrsfs}
\renewcommand\mathcal{\mathscr}

\newcommand{\M}{{\cal M}}

\newcommand{\G}{{\cal G}}

\newcommand{\C}{{\cal C}}

\newcommand{\maths}[1]{{\mathbb #1}}  

\newcommand{\RR}{\maths{R}}
\newcommand{\NN}{\maths{N}}

\newcommand{\SSS}{\maths{S}}

\newcommand{\ZZ}{\maths{Z}}

\newcounter{fig}

\newcommand{\dddp}{\partial_{\infty}^2\widetilde{\Sigma}}
\newcommand{\ddp}{\partial_{\infty}\widetilde{\Sigma}}

\newcommand{\Image}{\operatorname{Image}}

\newcommand{\Supp}{\operatorname{Supp}}

\newcommand{\CAT}{\operatorname{CAT}}

\newcommand{\revet}{(\widetilde{\Sigma},[\widetilde{q}])}

\newcommand{\srfce}{({\Sigma},{[q]})}

\newcommand{\revetm}{(\widetilde{\Sigma},\widetilde{m})}
\newcommand{\srfcem}{({\Sigma},m)}
\newcommand{\Id}{\operatorname{Id}}

\newcommand{\Lq}{\Lambda_{[q]}}

\newcommand{\Lqr}{\widetilde{\Lambda}_{[\widetilde{q}]}}
\newcommand{\Lm}{\Lambda_m}

\newcommand{\Lmr}{\widetilde{\Lambda}_{\widetilde{m}}}
\newcommand{\Gqr}{\G_{[\widetilde{q}]}}
\newcommand{\Gmr}{\G_{\widetilde{m}}}
\newcommand{\Gq}{\G_{[q]}}

\newcommand{\mutilde}{\widetilde{\mu}}
\newcommand{\grperevet}{\Gamma_{\widetilde{\Sigma}}}

\title{Measured geodesic laminations in Flatland.}
\author{Thomas~Morzadec}

\begin{document}
\maketitle
 
\textbf{Abstract:} Since their introduction by Thurston, measured geodesic laminations on hyperbolic surfaces occur in many contexts. 
In this survey (see \cite{Morzy1,Morzy2} for a complete exposition and proofs), we give a 
generalization of geodesic laminations on surfaces endowed with a half-translation structure (that is a singular flat surface with holonomy 
$\{\pm\Id\}$), called {\it flat laminations}, and we define transverse measures on flat laminations similar to
transverse measures on hyperbolic laminations, taking into account that the images of the leaves of a flat lamination are in general not pairwise disjoint. 
One aim is to construct a tool that could allow a fine description of the space of degenerations of half-translation structures on a surface. 
We define a topology on the set of measured flat laminations and a natural continuous 
projection of the space of measured 
flat laminations onto the space of
measured hyperbolic laminations, for any arbitrary half-translation structure and hyperbolic metric on a surface. We prove in particular that the space of measured flat
laminations is projectively compact. The main result of this survey is a classification theorem of (measured) flat laminations on a compact surface endowed with a 
half-translation structure. We also give an exposition of that every finite metric fat graph, outside four homeomorphisms classes, is the support of uncountably many
 measured flat laminations with uncountably 
many leaves none of which is eventually periodic, and that the space of measured flat laminations is separable and projectively compact.     
\footnote{  
AMS codes 30F30, 53C12, 53C22.}

%

\section{Introduction.}

The main aim of this paper is to survey a generalization of geodesic laminations on hyperbolic surfaces (see for instance \cite{Bonahon97}) to the surfaces endowed
with a half-translation structure, which is a flat metric with conical singular points and with holonomies in $\{\pm\Id\}$, that we will call {\it flat laminations},
as well as a definition of transverse measures on flat laminations. We refer to \cite{Morzy1,Morzy2} for complete proofs. Although the definitions are inspired of
measured geodesic laminations on hyperbolic surfaces, the extension is non trivial because of several surprising phenomena, notably since the images of the leaves of a flat 
lamination are in general not pairwise disjoint.
We will call {\it measured flat lamination} a flat lamination endowed with a transverse measure. We will define a sufficiently
fine
topology on the set of measured flat laminations. We will construct a (non injective) natural continuous projection of the space of measured flat laminations onto 
the space of
measured
hyperbolic laminations, for any choice of half-translation structure and (complete) hyperbolic metric on a surface.
We will completely describe the fibers not reduced to a point: they are the set of simple closed hyperbolic geodesics which are homotopic to a non trivial flat cylinder 
of the surface (endowed with the half-tranlation structure).

\medskip

The space of measured flat laminations allows to consider the 
measured flat laminations that are the limits of some sequences of periodic local geodesics, in the projectivized space of measured flat laminations. 
This in turn could yield a better understanding of the degenerations of half-translation structures on a surface, as initiated in \cite{DucLeiRaf10}. 
In particular, as spaces of measures are suitable for analysis tools (distributions as in \cite{Bonahon97}), this could allow a finer study of the boundary of the space of
half-translation structures that we will develop in a subsequent work.

Let $\Sigma$ be a compact, connected, orientable surface. In this survey, we assume that $\Sigma$ is without boundary,
and we refer to \cite{Morzy1,Morzy2} for the extension. 
A {\it half-translation structure} (or flat structure with conical singularities and holonomies in $\{\pm\Id\}$) on $\Sigma$ is the data consisting in a 
(possibly empty) discrete set of points $Z$ of $\Sigma$ and of a
Euclidean metric on $\Sigma-Z$ with conical {\it singular points} of angles of the form $k\pi$, with $k\in\NN$ and $k\geqslant 3$, at each point of $Z$, such
that the holonomy of every
piecewise $\C^1$ loop of $\Sigma-Z$ is contained in $\{\pm\Id\}$.

The surface $\Sigma$ endowed with a half-translation structure is a complete and locally $\CAT(0)$ metric space $(\Sigma,d)$. Let $p:(\widetilde{\Sigma},\widetilde{d})\to
(\Sigma,d)$ be a locally isometric universal cover. Two local geodesics ${\ell},{\ell}'$ of $(\Sigma,d)$, defined up to 
changing the
origins, are said to be 
{\it interlaced} if they have some lifts $\widetilde{\ell},\widetilde{\ell}'$ in $\widetilde{\Sigma}$ such that the image of $\widetilde{\ell}$ intersects both
complementary
components of $\widetilde{\ell}'(\RR)$ in $\widetilde{\Sigma}$, and conversely. A local geodesic is said to be {\it self-interlaced} if it is interlaced with itself. 
We endow the set of oriented, but non parametrized, local geodesics of $(\Sigma,d)$ with the quotient topology of the compact-open topology for the action by translations
on the parametrizations,
of $\RR$ on the parametrized local geodesics, which is called the {\it geodesic topology}. 

\bdefi A (geodesic) flat lamination on $(\Sigma,d)$ is a non empty set $\Lambda$ of complete local geodesics of $(\Sigma,d)$, defined up to changing 
origin, whose elements are called {\it leaves}, such that:

\begin{itemize}
 \item[$\bullet$]the leaves of $\Lambda$ are non self-interlaced and pairwise non interlaced;
 \item[$\bullet$] $\Lambda$ is invariant by changing the orientations of the leaves;
 \item[$\bullet$] $\Lambda$ is closed for the geodesic topology.
\end{itemize}

We will call {\it support} of $\Lambda$ the union of the images of the leaves of $\Lambda$. 
\edefi

Here are some examples of flat laminations on $\Sigma$. A {\it cylinder lamination} is a closed set of parallel leaves (and their opposites) whose images
are contained in a non degenerated flat cylinder (hence, 
these leaves are periodic). A {\it minimal flat lamination} is a flat lamination which is the closure,
for the geodesic topology, of the union of a leaf $\ell$ and its opposite $\ell^-$. It is of 
{\it recurrent type} if $\ell$ is regular
(i.e. does not meet any singular point) and $\ell$ is not periodic.
All the images of its regular leaves are then the images of some regular leaves of the vertical foliation
of a quadratic differential, and thus are dense in a {\it domain} of $\Sigma$, i.e. the closure 
of a connected open subset bounded by periodic local geodesics. And it is of {\it finite graph type} if the image of
$\ell$ is a finite graph, and if neither $\ell$ nor its 
opposite are periodic after a certain time. All the images of its leaves are then equal, and no leaf is periodic after a certain time. For example,
a minimal vertical foliation
of a holomorphic quadratic differential $q$ on a compact, connected Riemann surface is a minimal flat lamination of recurrent type, for the half-tranlation structure 
defined by $q$ (with a little work for the non regular leaves, see \cite[Lem.~4.11]{Morzy1}). We will construct minimal flat laminations of finite graph type in Theorem
\ref{thesecond'}.

We say that an end of a leaf (see Section \ref{rappelsconventions} for the definition) {\it terminates} in a minimal lamination or in a cylinder lamination
if a representative ray is equal to a ray of
a leaf of the lamination
or of a boundary component of the cylinder containing the support of the cylinder lamination, up to changing the origin. The main results of this paper are the following two theorems.

\btheo Every flat lamination on $\Sigma$ is a finite union of cylinder components, of minimal components (of recurrent type, finite graph type and periodic 
leaf travelled in both orientations)
and of isolated leaves (for the geodesic topology) both of whose ends terminate in a minimal component or a cylinder component.
\etheo

New phenomena appear in flat laminations compared with hyperbolic ones: the images of two leaves are generally not disjoint, the flat laminations are
not determined by their supports (uncountably many flat laminations can have the same support,
and contrarily to a hyperbolic lamination, a flat lamination may not be minimal whereas the image of each of its leaves is dense in the support of the lamination),
the cylinder components may have uncountably many
leaves.
Finally, there are three types, and no longer two, of minimal components of a flat lamination on a compact surface (periodic leaf travelled in both orientations,
minimal component of 
recurrent type or of finite graph type, see Theorem \ref{1Morzy} for a complete statement). Compared with hyperbolic laminations, the main difficulty to define
transverse 
measures on flat laminations is that the images of the leaves are not necessarly disjoint and that the support does not determine the lamination.
Hence, we no longer define
the transverse measure
as a family of measures on the images of the arcs transverse to the lamination, but as a family of measures on the sets of local geodesics that
intersect them transversally,
and we have to refine the notion of invariance by holonomy of these families of measures.

A {\it cyclic orientation} on a finite metric graph $X$ is the data of a  cyclic order (see \cite[§~2.3.1]{Wolf11} for the definition) on the set of germs of edges
starting at each vertex of the graph.

\btheo\label{thesecond'} Every cyclically oriented, connected, finite, metric graph $X$, without terminal point , is the support of uncountably many
uncountable minimal flat laminations without
eventually periodic leaf,  on a 
compact, connected, orientable surface endowed with a  half-translation structure, except if $X$ is homeomorphic to a circle, a dumbbell pair,
a flat height or a flat theta, by a homeomorphism preserving the cyclic
orientations (i.e  \begin{picture}(0,0)%
\includegraphics{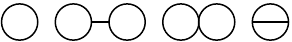}%
\end{picture}%
\setlength{\unitlength}{3771sp}%
\begingroup\makeatletter\ifx\SetFigFont\undefined%
\gdef\SetFigFont#1#2#3#4#5{%
  \reset@font\fontsize{#1}{#2pt}%
  \fontfamily{#3}\fontseries{#4}\fontshape{#5}%
  \selectfont}%
\fi\endgroup%
\begin{picture}(1461,198)(-8,830)
\end{picture}%
, where the orientations are given by the plan).
\etheo   

\medskip

Since $\Sigma$ is compact, if $m$ is a hyperbolic metric on $\Sigma$ and $\widetilde{m}$ is the pull back of $m$ on $\widetilde{\Sigma}$, there exists a unique $\grperevet$-equivariant homeomorphism between
the boundaries at infinity of $\revet$ and of $\revetm$, that allows to identify them. Let $\ddp$ denote the boundary at infinity of $\widetilde{\Sigma}$ and
$\dddp=\ddp\times\ddp-\{(x,x),\,x\in\ddp\}$. If $g$ is a geodesic of $\revet$ or of $\revetm$, let $E(g)=(g(-\infty),g(+\infty))$ be its ordered pair of points at infinity.  
Then,
for every $(x,y)\in\dddp$, there exists a unique geodesic $\widetilde{\lambda}$
of $\revetm$ such that $E(\widetilde{\lambda})=(x,y)$. However, it always exists a geodesic $\widetilde{\ell}$ of $\revet$ such that $E(\widetilde{\ell})=(x,y)$,
but it may happen that this geodesic is not unique, and then the set of such geodesics of $\revet$ is a set of parallel geodesics foliating a
flat strip,
and they project to periodic local geodesics foliating a flat cylinder. This point allows to define a non bijective correspondance between flat and hyperbolic measured laminations, from which we deduce 
many results of this survey and notably the two theorems above. We will also give an exposition of that the space of measured flat laminations is separable and projectively compact.

\medskip

In Section \ref{rappelsconventions}, we define the flat laminations and the transverse measures on flat laminations, and we endow the set of measured flat laminations with a topology.
In Section \ref{liensplatshyperboliques}, we define  a proper, surjective, continuous map from the space of measured flat laminations
to the space of measured hyperbolic laminations, for a complete hyperbolic metric, and we characterize its lack of injectivity. From this, we give a sketch of proof
of Theorem \ref{thesecond'}. In Section \ref{arbredual}, we define the tree associated to a measured flat lamination on a compact surface, and we construct
the covering group action on it.

\medskip\noindent
{\small{\it Acknoledgement: I want to thank Frederic Paulin for many advice and corrections that have deeply improved the redaction of this survey.}

\section{Definitions.}\label{rappelsconventions}

In this section \ref{rappelsconventions}, we give the definitions of half-translation structures on surfaces, of flat laminations and of transverse measures on flat laminations.

\subsection{Half-translation structures on a surface.}\label{structuresplates}

In the whole paper, we will use the definitions and notation of \cite{BriHae99} for a surface endowed with a distance
$(\Sigma,d)$: (locally) $\CAT(0)$, $\delta$-hyperbolic,... 
Notably, a {\it geodesic} (resp. a {\it local geodesic}) 
of $(\Sigma,d)$ is an isometric (resp. locally isometric) map $\ell:I\to \Sigma$, where $I$ is an interval of $\RR$. It will be called a 
  {\it segment}, a {\it ray} or a {\it  geodesic line} of $(\Sigma,d)$ if $I$ is respectively a compact interval, a closed half line 
  (generally $[0,+\infty[$)
  or $\RR$.
  If there is no precision, a {\it geodesic} is a geodesic line. A
 {\it germ of geodesic ray}, or simply a {\it germ}, is an equivalence class of locally geodesic rays for the equivalence relation  
 $r_1\sim_0 r_2$ if $r_1$ and $r_2$ coïncide on a non empty initial segment that is not reduced to a point. Similarly, the relation $r\sim_\infty r'$ if there exist $T,T'>0$
 such that $r(t+T)=r'(t+T')$ for all $t\geqslant 0$, is an equivalence relation on the set of subrays of a local geodesic.
An equivalence class for this equivalence relation is called an {\it end} (in the sense of Freudhental) of a local geodesic. A local geodesic has two ends. 
    \medskip

Let $\Sigma$ be a connected, orientable surface. In the whole article, we will assume that the boundary of $\Sigma$ is empty, for simplicity, but the results extend to
surfaces with non empty boundary (see \cite{Morzy1,Morzy2}).

\bdefi
A half-translation structure (or flat structure with conical singularities and holonomies in $\{\pm\Id\}$) on a surface $\Sigma$ is the data of a (possibly empty) discrete
subset $Z$ of $\Sigma$
and a Euclidean metric on $\Sigma-Z$ with conical singularity of angle $k_z\pi$ at each $z\in Z$, with $k_z\in\NN$ and $k_z\geqslant 3$,
such that the holonomy of every piecewise $\C^1$ loop in 
 $\Sigma-Z$ is contained in $\{\pm\Id\}$.
\edefi

Two vectors $v_1$ and $v_2$ tangent to 
$\Sigma$ have the 
{\it same direction} if $v_2$ is the image of $\pm v_1$ by holonomy along a piecewise $\C^1$ path of $\Sigma-Z$ between the basepoints of 
$v_1$ and $v_2$. 
This definition does not depend on the choice of a path, since the holonomy of every loop is contained in $\{\pm\Id\}$. Hence, there is a notion of direction on a surface
endowed
with a half-translation structure, but there is no "vertical direction", as defined by a quadratic differential. A piecewise $\C^1$ path or union of paths is said to have
{\it constant direction}, if all its tangent vectors, at the points in $\Sigma-Z$, have the same direction.

\medskip

 We will denote by $[q]$ a half-translation structure on $\Sigma$, with $q$ a holomorphic quadratic differential (see \cite[§~2.5]{Morzy1} for an explanation of
 the notation).
 A half-translation structure defines a geodesic distance $d$ that is locally $\CAT(0)$. We will call {\it local flat geodesics} the local geodesics of a half-translation 
 structure. A continuous map $\ell:\RR\to\Sigma$ is a local flat geodesic if and only if it satisfies (see \cite[Th.~5.4~p.24]{Strebel84} and
 \cite[Th.~8.1~p.~35]{Strebel84}): for every $t\in\RR$,

\medskip
\noindent
$\bullet$~ if $\ell(t)$ does not belong to $Z$, there exists a neighborhoud $V$ of $t$ in $\RR$ such that $\ell_{|V}$ is an Euclidean segment
(hence, $\ell_{|V}$ has constant direction);
 
\noindent
$\bullet$~ if $\ell(t)$ belongs to $Z$, then the two angles defined by the germs of $\ell([t,t+\varepsilon[)$ and $\ell(]t-\varepsilon,t])$, 
with $\varepsilon>0$ small enough, measured in both connected components of $U-\ell(]t-\varepsilon,t+\varepsilon[)$, with $U$  a small enough neighborhoud of
$\ell(t)$, are at least $\pi$.
\begin{center}
\begin{picture}(0,0)%
\includegraphics{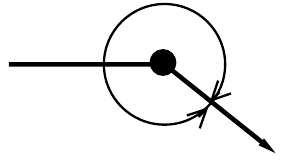}%
\end{picture}%
\setlength{\unitlength}{5801sp}%
\begingroup\makeatletter\ifx\SetFigFont\undefined%
\gdef\SetFigFont#1#2#3#4#5{%
  \reset@font\fontsize{#1}{#2pt}%
  \fontfamily{#3}\fontseries{#4}\fontshape{#5}%
  \selectfont}%
\fi\endgroup%
\begin{picture}(922,517)(12226,-1181)
\put(12241,-1051){\makebox(0,0)[lb]{\smash{{\SetFigFont{14}{16.8}{\familydefault}{\mddefault}{\updefault}{\color[rgb]{0,0,0}$\pi\leqslant$}%
}}}}
\put(13006,-871){\makebox(0,0)[lb]{\smash{{\SetFigFont{14}{16.8}{\familydefault}{\mddefault}{\updefault}{\color[rgb]{0,0,0}$\geqslant\pi$}%
}}}}
\end{picture}%

\end{center}

\subsection{Geodesic laminations on surfaces endowed with a half-translation structure and with a (complete) hyperbolic metric.}

Let $\Sigma$ be a connected, orientable surface (without boundary). Let $[q]$ be a half-translation structure and let $m$ be a complete hyperbolic metric on $\Sigma$. Let  
$p:\widetilde{\Sigma}\to\Sigma$ be a  universal cover of covering group $\Gamma_{\widetilde{\Sigma}}$, let $[\widetilde{q}]$ be the unique half-translation structure 
and let $\widetilde{m}$ be the unique hyperbolic metric  
on $\widetilde{\Sigma}$ such that $p:\revet\to\srfce$ and $p:\revetm\to\srfcem$ are locally isometric.

\medskip

We call  
{\it geodesic topology} the compact-open topology on the set $\G_d$ of parametrized local geodesics for the distance $d$ defined by $[q]$ or by $m$,
   or the quotient topology of the compact-open topology by the action by translations on the parametrizations, of $\RR$, on the set $[\G_d]$ of local geodesics defined
   up to changing origin. They will be the only topologies taken into account on the spaces of local geodesics. The quotient projection from $\G_d$ to $[\G_d]$  will be
   denoted by $g\mapsto[g]$, and if $F\subset\G_d$, we will denote by $[F]$ its image in $[\G_d]$. If $g:\RR\to (\Sigma,d)$ is an element of $\G_d$ or of $[\G_d]$,
   we denote its {\it opposite} by $g^-(t):t\mapsto g(-t)$. Finally, if $F$ is a set of elements of $\G_d$ or of $[\G_d]$, we call {\it support of $F$} the union
   of the images of the elements of $F$, denoted by $\Supp(F)$. 
   
   \blemm\cite[Coro.~2.4]{Morzy1} If $F$ is a closed set of $\G_d$ or of $[\G_d]$ (for the geodesic topology), then $\Supp(F)$ is closed in $\Sigma$.
   \elemm
    
    The opposite is true on a hyperbolic surface, but it may be false on a surface endowed with a half-translation structure. 
    We will see (Theorem \ref{thesecond}) that it exists some surfaces endowed with a half-tranlation structure, having a flat local geodesic $\ell$ whose
    image is a
    finite graph (and hence is closed), such that
    the closure of $\ell$ for the geodesic topology contains uncountably many geodesics.
    
    \medskip

 In \cite[§~2.3]{Morzy1}, we have given a very global definition of {\it interlaced local geodesics}, in a locally $\CAT(0)$, complete, connected metric space, 
 whose boundary at infinity of a universal cover is endowed with a (total) cyclic order. Here, we only recall the specific definition in the case of
 connected, orientable surfaces without boundary, endowed with a complete locally $\CAT(0)$ metric $d$.  
  Since the  universal cover $(\widetilde{\Sigma},\widetilde{d})$ is  $\CAT(0)$, the intersection of the images of two geodesics of $(\widetilde{\Sigma},\widetilde{d})$ is connected, possibly 
  empty. If $\widetilde{\ell}$ is a geodesic of $(\widetilde{\Sigma},\widetilde{d})$, defined up to changing the origin, then $\widetilde{\Sigma}-\widetilde{\ell}(\RR)$
  has two connected components. Two complete geodesics $\widetilde{\ell},\widetilde{\ell'}$ of $(\widetilde{\Sigma},\widetilde{d})$, defined up to changing the origins, are
{\it interlaced} if 
$\widetilde{\ell}$ meets both connected components of $\widetilde{\Sigma}-\widetilde{\ell}'(\RR)$, or the opposite (which is equivalent). Two local geodesics of  
$(\Sigma,d)$ are {\it interlaced} if they admit some lifts in $(\widetilde{\Sigma},\widetilde{d})$ which are interlaced, and a local geodesic is 
{\it self-interlaced} if it is interlaced with itself.

If $d$ is the distance defined by $m$, two local geodesics are non interlaced if and only if their images are disjoint and a local geodesic is non self-interlaced if and only if 
it is simple, but this equivalence is not true in the case of $[q]$. We refer to \cite[§~3.1]{Morzy1} for a characterization of 
the local geodesics for $[q]$ that are non interlaced.

\bdefi\label{laminationplate1}
A geodesic lamination (or simply a lamination) of $(\Sigma,d)$, with $d$ the distance defined by $m$ or $[q]$, is a non empty set  
$\Lambda$ of (complete) local geodesics of $(\Sigma,d)$, defined up to changing origin, whose elements are called leaves, such that: 

\medskip

\noindent
$\bullet$~ leaves are non self-interlaced;

\noindent
$\bullet$~ leaves are pairwise non interlaced;

\noindent
$\bullet$~ if $\ell$ belongs to $\Lambda$ then so do $\ell^-$, with $\ell^{-}(t):t\mapsto\ell(-t)$;

\noindent
$\bullet$~ $\Lambda$ is closed for the geodesic topology.
\edefi

We say that $\Lambda$ is a {\it  flat lamination} if $d$ is defined by $[q]$ and that $\Lambda$ is a {\it hyperbolic lamination} if $d$ is defined by $m$.
Usually, a hyperbolic lamination is defined as
a non empty closed subset of $\Sigma$, which is a disjoint union of images of simple hyperbolic local geodesics.
The definitions are equivalent in the case of hyperbolic laminations but not in the case of flat laminations. For example, it may happen that a flat  local geodesic is not simple whereas it is 
non self-interlaced (see \cite[§~4.1]{Morzy1}).

\subsection{Measured flat laminations.}

 Let $\srfce$ be a connected, orientable surface (without boundary), endowed with a half-translation structure. An {\it arc} is a piecewise $\C^1$ map
 $\alpha:[0,1]\to\Sigma$ which is a homeomorphism onto its image. Let $\Lambda$ be a flat lamination of $\srfce$. An arc $\alpha$ is {\it transverse}
 to a leaf or to a segment of leaf $\ell$ of $\Lambda$ if

\noindent
$\bullet$~$\alpha$ is transverse to $\ell$ in the complementary of the singular points of $[q]$ and of the singular points of $\alpha$;

\medskip
\noindent
$\bullet$~ for every singular point $x$ of $[q]$ or of $\alpha$ in $\Image(\ell)\cap\alpha(]0,1[)$, there exists a neighborhoud
$U$ of $x$ that is a topological disk, and a segment $S$ of $\ell$ such that $U-\Image(S)\cap U$ has two connected components and the connected components of
$U\cap(\alpha([0,1])-\{x\})$ are contained in different components of $U-\Image(S)\cap U$;
\begin{center}
\begin{picture}(0,0)%
\includegraphics{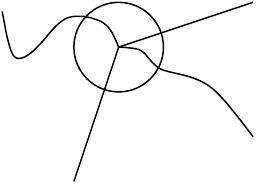}%
\end{picture}%
\setlength{\unitlength}{3771sp}%
\begingroup\makeatletter\ifx\SetFigFont\undefined%
\gdef\SetFigFont#1#2#3#4#5{%
  \reset@font\fontsize{#1}{#2pt}%
  \fontfamily{#3}\fontseries{#4}\fontshape{#5}%
  \selectfont}%
\fi\endgroup%
\begin{picture}(1284,924)(1204,-1873)
\put(1801,-1096){\makebox(0,0)[lb]{\smash{{\SetFigFont{11}{13.2}{\familydefault}{\mddefault}{\updefault}{\color[rgb]{0,0,0}$x$}%
}}}}
\put(1486,-1501){\makebox(0,0)[lb]{\smash{{\SetFigFont{11}{13.2}{\familydefault}{\mddefault}{\updefault}{\color[rgb]{0,0,0}$U$}%
}}}}
\put(2116,-1231){\makebox(0,0)[lb]{\smash{{\SetFigFont{11}{13.2}{\familydefault}{\mddefault}{\updefault}{\color[rgb]{0,0,0}$\ell$}%
}}}}
\put(2071,-1501){\makebox(0,0)[lb]{\smash{{\SetFigFont{11}{13.2}{\familydefault}{\mddefault}{\updefault}{\color[rgb]{0,0,0}$\alpha$}%
}}}}
\end{picture}%

\end{center}

\noindent
$\bullet$~ $\alpha$ is tangent to $\ell$  neither in $0$  nor in $1$.

\medskip

An arc $\alpha$ is {\it transverse to a set $F$ of leaves} or of segments of leaves of  $\Lambda$ if it is transverse to every element of $F$.
In particular, an arc is  transverse to $\Lambda$ if it is transverse to every leaf of $\Lambda$.

If $\alpha:[0,1]\to\Sigma$ is an arc of  
 $\Sigma$, we denote by $G(\alpha)$ the compact set of elements of $\Gq$ (the set of parametrized local geodesics of $\srfce$) which are transverse to $\alpha$ and whose origins 
 belong to 
 $\alpha([0,1])$. By definition, if $\alpha'([0,1])\subseteq\alpha([0,1])$, then $G(\alpha')\subseteq G(\alpha)$. Let $F_1\subseteq\Gq$ be such that $[F_1]\subseteq\Lambda$ and let $\alpha_1$ and $\alpha_2$ be two disjoint  arcs transverse to $F_1$, such that
$F_1\subseteq G(\alpha_1)$ and every element of $F_1$ intersects $\alpha_2([0,1])$ at a positive time. For every $g_1\in F_1$, we define $t_{g_1}=
\min\{t>0\;:\;g_1(t)\in \alpha_2([0,1])\}$, and $F_2$ the subset of the elements $g_2\in G(\alpha_2)$ such that there exists $g_1\in F_1$ with 
$g_2(t)=g_1(t+t_{g_1})$ for all $t\in\RR$. A {\it holonomy} $h:F_1\to F_2$ of $\Lambda$ is a homeomorphism between $F_1$ and $F_2$ defined by 
$h(g_1)=g_2:t\mapsto g_1(t+t_{g_1})$ such that there exists a homotopy $H:[0,1]\times[0,1]\to\Sigma$ between $\alpha_1$ and $\alpha_2$ such that:

\noindent
$\bullet$~ for every $t\in\,[0,1]$,  the map $s\mapsto H(s,t)$ is an arc transverse to every segment of leaf $g_{1|[0,\,t_{g_1}]}$, with $g_1\in
F_1$;

\noindent
$\bullet$~ for every $\ell\in F_1$, there exists $s_\ell\in[0,1]$ such that $t\mapsto H(s_\ell,t)$ is a segment of $\ell$ (up to changing the parametrization);

\noindent
$\bullet$~ the intersections $H([0,1]\times]0,1[)\cap \alpha_i([0,1])$ with $i=1,2$ are empty.

Contrarily to the case of measured foliations, if the images of the geodesics are not pairwise disjoint, the map $H$ may not be injective.

\bdefi\label{defmesuretransverse} 
A transverse measure on $\Lambda$ is a family $\mu=(\mu_\alpha)_\alpha$ of Radon measures $\mu_\alpha$ defined on $G(\alpha)$, for every  arc $\alpha$ transverse to $\Lambda$, such that:

\noindent
$(1)$~ the support of $\mu_\alpha$ is the set $\{\ell\in G(\alpha)\;:\;[\ell]\in\Lambda\}$;

\noindent
$(2)$~ if $h:F_1\to F_2$ is a holonomy of $\Lambda$, where $\alpha_1,\alpha_2$ are two disjoint  arcs transverse to $F_1$ and $F_1\subset G(\alpha_1)$ and 
$F_2\subset G(\alpha_2)$ are some Borel sets, then $h_*(\mu_{\alpha_1|F_1})=\mu_{\alpha_2|F_2}$;

\noindent
$(3)$~ $\mu_\alpha$ is $\iota$-invariant, with $\iota(\ell)=\ell^-:t\mapsto\ell(-t)$;

\noindent
$(4)$~ if $\alpha'([0,1])\subseteq\alpha([0,1])$, then $\mu_{\alpha|G(\alpha')}=\mu_{\alpha'}$. 
\edefi

We will denote by $(\Lambda,\mu)$ a flat lamination endowed with a transverse measure, that we will call a {\it measured flat lamination}, and we will denote by 
$\mathcal{M}\mathcal{L}_p(\Sigma)$ the set of measured flat laminations on $\Sigma$.
We endow $\mathcal{M}\mathcal{L}_p(\Sigma)$ with the topology such that a sequence $(\Lambda_n,\mu_n)_{n\in\NN}$ converges to $(\Lambda,\mu)$ if and only if for every arc 
$\alpha$, if $\alpha$ is transverse to $\Lambda$, then $\alpha$ is transverse to $\Lambda_n$ for $n$ large enough and $\mu_{n,\alpha}\overset{*}{\rightharpoonup}\mu_\alpha$ 
in the space of Radon measures on $G(\alpha)$.

A leaf $\ell$ of $\Lambda$ is {\it positively recurrent} if there exists an arc $\alpha$ transverse to ${\ell}$ such that $\ell$ intersects $\alpha([0,1])$ at an infinite 
number of positive times.

\blemm\cite[Lem.~7]{Morzy2}
If $\Lambda$ is endowed with a transverse measure $\mu$, then the only leaves of $\Lambda$ which are isolated and positively recurrent are the periodic leaves. 
\elemm

\section{Link between measured flat laminations and measured hyperbolic laminations.}\label{liensplatshyperboliques}

 In this section \ref{liensplatshyperboliques}, we denote by $\srfce$ a compact, connected, orientable surface (without boundary),
 endowed with a half-translation structure, and by
 $p:\revet\to\srfce$ a locally isometric universal cover whith covering group $\Gamma_{\widetilde{\Sigma}}$. We assume that $\chi(\Sigma)<0$,
 and we denote by $m$ a hyperbolic metric on $\Sigma$ and by $\widetilde{m}$ the pull back of $m$ on $\widetilde{\Sigma}$.
 
 Since $\Sigma$ is compact, there exists a unique $\grperevet$-equivariant homeomorphism between
the boundaries at infinity of $\revet$ and of $\revetm$, that allows to identify them. Let $\ddp$ denote the boundary at infinity of $\widetilde{\Sigma}$ and
$\dddp=\ddp\times\ddp-\{(x,x),\,x\in\ddp\}$. If $g$ is a geodesic of $\revet$ or of $\revetm$, lets $E(g)=(g(-\infty),g(+\infty))$ be its ordered pair of points at infinity. 
 
 If $\Lambda$ is a geodesic lamination of $\srfce$ (resp. $\srfcem$), and if $\widetilde{\Lambda}$ is the set of lifts of the leaves of $\Lambda$ 
 in $\widetilde{\Sigma}$, then $\widetilde{\Lambda}$ is a geodesic lamination of $\revet$ (resp. $\revetm$) which is $\Gamma_{\widetilde{\Sigma}}$-equivariant. 
 Conversely, if $\widetilde{\Lambda}$ is a $\Gamma_{\widetilde{\Sigma}}$-equivariant lamination of $\revet$ (resp. $\revetm$), then it projects to a lamination of
 $\srfce$ (resp. $\srfcem$), and the set $\mathcal{L}_{[q]}(\Sigma)$ (resp. $\mathcal{L}_h(\Sigma)$) of geodesic laminations of $\srfce$ (resp. $\srfcem$) is homeomorphic to
 the set $\mathcal{L}_{\grperevet}([\Gqr])$ (resp. $\mathcal{L}_{\grperevet}([\Gmr])$) of geodesic laminations of $\revet$ (resp. $\revetm$) which are 
 $\grperevet$-equivariant. Similarly,
 the spaces $\mathcal{M}\mathcal{L}_{[q]}(\Sigma)$ (resp. $\mathcal{M}\mathcal{L}_h(\Sigma)$) of measured flat (resp. hyperbolic) laminations on $\Sigma$, endowed respectivelly
 with the topology defined above and with the topology defined in \cite[p.~19]{Bonahon97}, are homeomorphic to the spaces of Radon measures  on $[\Gqr]$ (resp.
 $[\Gmr]$) which are $\Gamma_{\widetilde{\Sigma}}$ et $\iota$-invariant (with $\iota(\widetilde{g})=\widetilde{g}^-:t\mapsto\widetilde{g}(-t)$),
 whose supports are flat (resp. hyperbolic) laminations, denoted by $\M_{\grperevet}([\Gqr])$ (resp. $\M_{\grperevet}([\Gmr])$) (see
 \cite[Prop.~17~p.~154]{Bonahon88} and \cite[§4]{Morzy2}). 
 
 If $(x,y)\in\dddp$ there exists a unique geodesic $\widetilde{\lambda}$ of $\revetm$ such that $E(\widetilde{\lambda})=(x,y)$. It is not the same in $\revet$. There always
 exists a geodesic $\widetilde{\ell}$ of $\revet$ such that $E(\widetilde{\ell})=(x,y)$, but it may not be unique. 
 We recall that two geodesics $c$ and $c'$ are at \textit{finite Hausdorff distance} if there exists $K>0$ such that
$d(c(t),c'(t))\leqslant K$ for all $t\in\RR$.

\blemm\cite[Th.~2.(c)]{MS85}\label{geodasymptotes}
We recall that $\Sigma$ is compact. Let $\widetilde{\ell}_1$ and $\widetilde{\ell}_2$ be two geodesics of $\revet$, such that $\widetilde{\ell}_1$ and $\widetilde{\ell}_2$
are at finite Hausdorff distance and such that there exists $\delta>0$ with $d(\widetilde{\ell}_1(t),\widetilde{\ell}_2(\RR))\geqslant\delta$ for all $t\in\RR$.
Then, the convex hull of 
$\widetilde{\ell}_1(\RR)\cup \widetilde{\ell}_2(\RR)$ is isometric to a flat strip $\RR\times[0,D]\subset \RR^2$, with $D\geqslant 0$, and 
the projections $p\circ\widetilde{\ell}_1$ and $p\circ\widetilde{\ell}_2$ of $\widetilde{\ell}_1$ and $\widetilde{\ell}_2$ on $\Sigma$ are periodic local geodesics
which are freely homotopic to the boundary of a flat cylinder.
\elemm
 
Hence, if $\varphi:[\Gqr]\to[\Gmr]$ is the map associating to the geodesic $\widetilde{\ell}\in[\Gqr]$ the unique geodesic $\varphi(\widetilde{\ell})\in[\Gmr]$
such that $E(\varphi(\widetilde{\ell}))=E(\widetilde{\ell})$
  (see \cite[§4.2]{Morzy1}), then $\varphi$ is surjective and continuous, and a closed subset $F$ of $[\Gqr]$ is a flat lamination if and only if $\varphi(F)$ is a 
  hyperbolic
  lamination. However, the map $\varphi$ is not injective. By definition, two different geodesics $\widetilde{\ell},\widetilde{\ell}'$ of $\revet$ have the same
  image by
  $\varphi$ if and only if they have the same ordered pair of points at infinity. According to Lemma \ref{geodasymptotes}, their images are thus parallel and contained in a 
  maximal flat strip of $\revet$, and their projections in $\srfce$  are freely homotopic 
  periodic local geodesics, hence their images are contained in a maximal flat cylinder. The points at infinity of
  $\widetilde{\ell}$, $\widetilde{\ell}'$ and $\varphi(\widetilde{\ell})=\varphi(\widetilde{\ell}')$  are hence the attractive and repulsive fixed points
  of an element of the covering group, 
  and the projection of $\varphi(\widetilde{\ell})$ in $\srfcem$ is a simple closed local geodesic. Hence, the restriction of $\varphi$ to the set of geodesics of  
  $\revet$ whose projections are not periodic is injective and its image is the set of geodesics of $\revetm$ whose projections are not closed local geodesics,
  and the preimage of a geodesic $\widetilde{\lambda}$ of $\revetm$ whose projection  is a simple closed local geodesic, is the set of geodesics of $\revet$ having the same
  ordered pair of points at infinity than $\widetilde{\lambda}$, which may be a unique geodesic, or may be a set of parallel geodesics whose images are contained in a flat
  strip, and whose projections are some periodic local geodesics freely homotopic to $p\circ\widetilde{\lambda}$.

  Since the sets $\mathcal{L}\srfce$ and $\mathcal{L}\srfcem$ of flat and hyperbolic laminations on $\Sigma$ are 
  respectivelly homeomorphic to
 the sets $\mathcal{L}_{\grperevet}\revet$ and $\mathcal{L}_{\grperevet}\revetm$ of $\grperevet$-invariant flat and hyperbolic laminations on $\widetilde{\Sigma}$,
 the map $\varphi$ defines a continous, surjective map $\psi:\mathcal{L}\srfce\to\mathcal{L}\srfcem$. It implies the structure theorem of flat laminations.

We recall the two main results of \cite{Morzy1} about flat laminations. In Theorem \ref{1Morzy}, $\Lambda$ is a flat lamination. A {\it cylinder component} is a maximal set of
leaves of $\Lambda$ whose images are contained in a non degenerated flat cylinder (hence, 
these leaves are periodic).
A {\it minimal component} is a sublamination which is the closure, for the geodesic topology, of a leaf $\ell$ and its opposite $\ell^-$. The minimal component is of 
{\it recurrent type} if $\ell$ is regular (i.e. does not meet any singular point) and is not periodic. All the images of its regular leaves are then the images of some
regular leaves of the vertical foliation
of a quadratic differential, and thus are dense in a {\it domain} of $\Sigma$, i.e. the closure 
of a connected open subset bounded by periodic local geodesics. The minimal component is of {\it finite graph type} if the image of
$\ell$ is a finite graph, and if neither $\ell$ nor its 
opposite are eventually periodic. All the images of its leaves are then equal, and no leaf is eventually periodic. 
An end of a local geodesic {\it terminates} in a minimal component or in a cylinder component if there exists a ray in the equivalence class of this end which is the ray of
a leaf of the minimal component or a ray of a boundary component of the corresponding flat cylinder.

\begin{picture}(0,0)%
\includegraphics{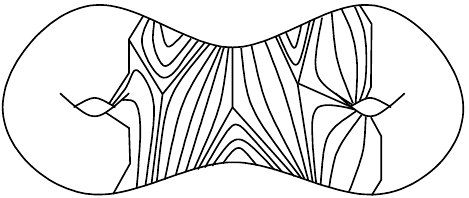}%
\end{picture}%
\setlength{\unitlength}{4144sp}%
\begingroup\makeatletter\ifx\SetFigFont\undefined%
\gdef\SetFigFont#1#2#3#4#5{%
  \reset@font\fontsize{#1}{#2pt}%
  \fontfamily{#3}\fontseries{#4}\fontshape{#5}%
  \selectfont}%
\fi\endgroup%
\begin{picture}(2130,896)(2449,-5534)
\end{picture}%
\;\begin{picture}(0,0)%
\includegraphics{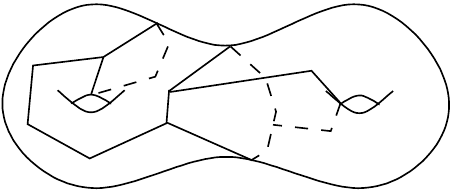}%
\end{picture}%
\setlength{\unitlength}{4144sp}%
\begingroup\makeatletter\ifx\SetFigFont\undefined%
\gdef\SetFigFont#1#2#3#4#5{%
  \reset@font\fontsize{#1}{#2pt}%
  \fontfamily{#3}\fontseries{#4}\fontshape{#5}%
  \selectfont}%
\fi\endgroup%
\begin{picture}(2066,865)(2472,-7830)
\end{picture}%
\;\begin{picture}(0,0)%
\includegraphics{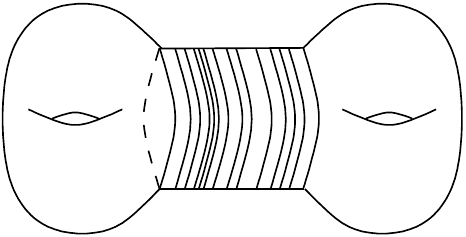}%
\end{picture}%
\setlength{\unitlength}{4144sp}%
\begingroup\makeatletter\ifx\SetFigFont\undefined%
\gdef\SetFigFont#1#2#3#4#5{%
  \reset@font\fontsize{#1}{#2pt}%
  \fontfamily{#3}\fontseries{#4}\fontshape{#5}%
  \selectfont}%
\fi\endgroup%
\begin{picture}(2124,1075)(3954,-5634)
\end{picture}%

\btheo\label{1Morzy}\cite[§~6]{Morzy1} Let $\Lambda$ be a flat lamination on a compact, connected, orientable surface endowed with a half-translation structure.
 Then $\Lambda$ is a finite union of cylinder components, of minimal components (of recurrent type, finite graph type and periodic leaf travelled in both orientations)
and of isolated leaves (for the geodesic topology) both of whose ends terminate in a minimal component or a cylinder component.
\etheo

This theorem comes from the structure theorem of hyperbolic laminations and from a few observations (see \cite[Coro.~4.6, Lem.~4.10, Lem.~4.11]{Morzy1}): 
the preimage of a minimal lamination of $\srfcem$, which is not
a closed leaf, is a minimal lamination of $\srfce$ which either contains a leaf whose image is not compact and then the lamination is of recurrent type, or contains a leaf
whose
image is compact, and then the lamination is of finite graph type. The preimage by $\psi$ of a closed leaf of $\srfcem$ is either a periodic leaf or a cylinder lamination. Finally, the preimage 
of an isolated local geodesic of $\srfcem$ whose an end spirales to a minimal hyperbolic lamination is an isolated local geodesic of $\srfce$ whose the corresponding end terminates
in the preimage of the hyperbolic lamination, since the flat geodesics can change direction at most at singular points (see \cite[§4.4]{Morzy1}).

\begin{center}
  \begin{picture}(0,0)%
\includegraphics{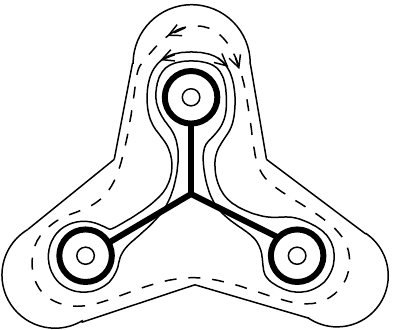}%
\end{picture}%
\setlength{\unitlength}{3771sp}%
\begingroup\makeatletter\ifx\SetFigFont\undefined%
\gdef\SetFigFont#1#2#3#4#5{%
  \reset@font\fontsize{#1}{#2pt}%
  \fontfamily{#3}\fontseries{#4}\fontshape{#5}%
  \selectfont}%
\fi\endgroup%
\begin{picture}(1959,1644)(4755,-2949)
\end{picture}%

 \end{center}
 
\rem If $\Lambda$ is a minimal flat lamination, then the image of every ray of every leaf of $\Lambda$ is dense in $\Supp(\Lambda)$.
However, contrary to the case of hyperbolic lamination, the opposite is false. For example, let consider the surface
endowed with a half-translation structure above, whose singular points are of angle $3\,\pi$ and located at the vertices of the graph. Then, the union of both local flat geodesics
in the free homotopy classes of the two closed curves drawn on the picture (and their opposites) is a flat lamination such that the image of every ray of every leaf is
equal to the graph, but the lamination has two minimal components. 

\medskip

We now show that every
finite, connected graph, without terminal point , outside four homeomorphism classes, is the support of uncountably many geodesic
laminations with uncountably many leaves none of which is eventually periodic.
 We use the definitions and notation about graphs of \cite{Serre83}. Every graph is identified with its topological realization, and it is endowed with a geodesic
distance such that each half-edge is isometric to a compact interval of $\RR$. In this paper, we only consider finite, connected graph, without terminal point.
 A {\it cyclic orientation} on a finite metric graph $X$ is the data of cyclic orders on the sets of half-edges
 issued from each vertex of $X$ (a cyclic order on the set of half-edges issued
from a vertex is globally an embedding of this set in the oriented unit circle $\SSS^1$, see \cite[§~2.3.1]{Wolf11} for the exact definition). A {\it cyclically oriented}
graph is a graph endowed with a cyclic orientation. In \cite[§5]{Morzy1},
we have recalled the classical result that every cyclically oriented, finite, connected graph $X$, without terminal point, can be isometrically embedded in a compact, 
connected 
surface $\Sigma(X)$ endowed 
with a half-translation structure, made of flat cylinders locally isometrically glued on the graph $X$, and that retracts by strong 
defomation on $X$ (that is why such a graph is often called a {\it fat graph}). 

\begin{center}
 \begin{picture}(0,0)%
\includegraphics{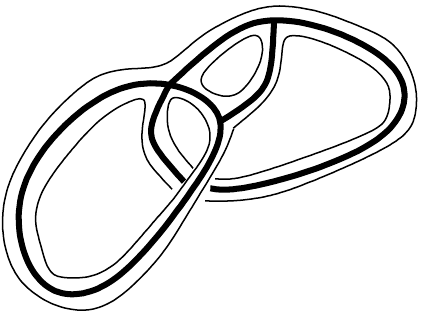}%
\end{picture}%
\setlength{\unitlength}{3812sp}%
\begingroup\makeatletter\ifx\SetFigFont\undefined%
\gdef\SetFigFont#1#2#3#4#5{%
  \reset@font\fontsize{#1}{#2pt}%
  \fontfamily{#3}\fontseries{#4}\fontshape{#5}%
  \selectfont}%
\fi\endgroup%
\begin{picture}(2082,1536)(4666,840)
\end{picture}%

\end{center}

 The Euler characteristic of this surface is negative, except if $X$ is a circle (which
embbeds in a flat cylinder), hence it can be endowed with a hyperbolic metric. If $\Sigma(X)$ is not a pair of pant, we have seen (see \cite[Lem.~4.2]{Morzy1}) 
that there exists a minimal hyperbolic
lamination $\Lambda$ on $\Sigma(X)$ (endowed with a complete hyperbolic metric with totally geodesic boundary) such that the completions of the connected components of
$\Sigma(X)-
\Supp(\Lambda)$ are ideal triangles or ideal monogons minus a disk (actually, the proof of the lemma implies that there exists uncountably many such hyperbolic lamination,
with a fixed hyperbolic metric). Then, we show that the preimage by $\psi$ of such a lamination is a minimal flat lamination, which is not a 
periodic leaf travelled in both orientations, whose support is the graph $X$.  Since $\Sigma(X)$ is a pair of pant if and only if $X$ is homeomorphic to a circle, a
dumbbell pair, 
a flat height or a theta by a homeomorphism preserving the cyclic orientations, we get the following theorem.

\begin{center}
 \begin{picture}(0,0)%
\includegraphics{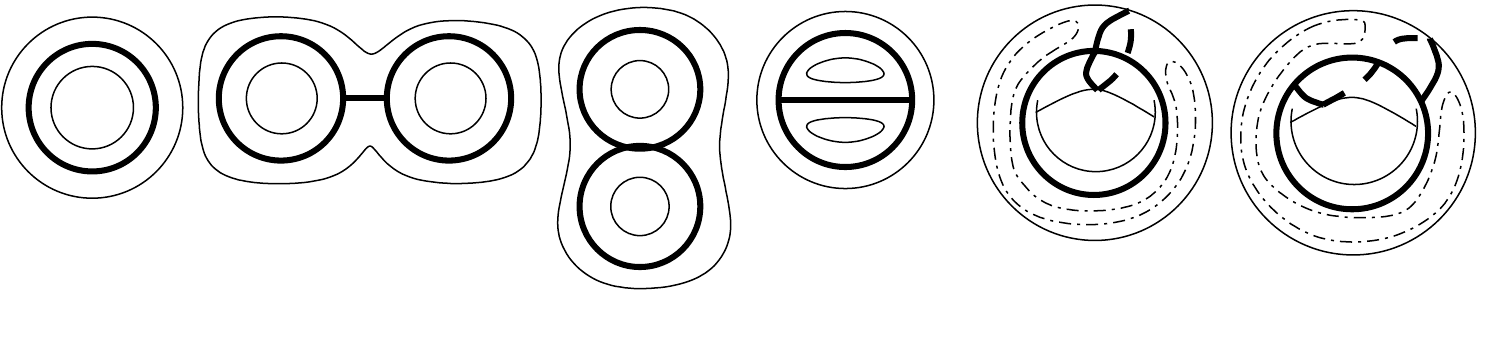}%
\end{picture}%
\setlength{\unitlength}{3812sp}%
\begingroup\makeatletter\ifx\SetFigFont\undefined%
\gdef\SetFigFont#1#2#3#4#5{%
  \reset@font\fontsize{#1}{#2pt}%
  \fontfamily{#3}\fontseries{#4}\fontshape{#5}%
  \selectfont}%
\fi\endgroup%
\begin{picture}(7386,1657)(683,1215)
\put(946,1739){\makebox(0,0)[lb]{\smash{{\SetFigFont{11}{13.2}{\familydefault}{\mddefault}{\updefault}{\color[rgb]{0,0,0}Circle}%
}}}}
\put(2206,1784){\makebox(0,0)[lb]{\smash{{\SetFigFont{11}{13.2}{\familydefault}{\mddefault}{\updefault}{\color[rgb]{0,0,0}Dumbbell}%
}}}}
\put(3556,1289){\makebox(0,0)[lb]{\smash{{\SetFigFont{11}{13.2}{\familydefault}{\mddefault}{\updefault}{\color[rgb]{0,0,0}Flat height}%
}}}}
\put(4546,1739){\makebox(0,0)[lb]{\smash{{\SetFigFont{11}{13.2}{\familydefault}{\mddefault}{\updefault}{\color[rgb]{0,0,0}Flat theta}%
}}}}
\put(5626,1514){\makebox(0,0)[lb]{\smash{{\SetFigFont{11}{13.2}{\familydefault}{\mddefault}{\updefault}{\color[rgb]{0,0,0}Torical height}%
}}}}
\put(6931,1469){\makebox(0,0)[lb]{\smash{{\SetFigFont{11}{13.2}{\familydefault}{\mddefault}{\updefault}{\color[rgb]{0,0,0}Torical theta}%
}}}}
\end{picture}%

\end{center}

\btheo\label{thesecond}\cite[Th.~5.2]{Morzy1} Every cyclically oriented, connected, finite, metric graph $X$, without extremal point, may be the support of uncountably many uncountable
minimal flat laminations
with no eventually periodic leaf,  on a 
compact and connected surface endowed with a  half-translation structure, except if $X$ is homeomorphic to a circle, a dumbbell pair, a flat height or a flat theta, 
by a homeomorphism preserving the cyclic
orientations (i.e \begin{picture}(0,0)%
\includegraphics{figure22.pdf}%
\end{picture}%
\setlength{\unitlength}{3771sp}%
\begingroup\makeatletter\ifx\SetFigFont\undefined%
\gdef\SetFigFont#1#2#3#4#5{%
  \reset@font\fontsize{#1}{#2pt}%
  \fontfamily{#3}\fontseries{#4}\fontshape{#5}%
  \selectfont}%
\fi\endgroup%
\begin{picture}(1461,198)(-8,830)
\end{picture}%
, where the orientations are given by the plan).
\etheo

 We now focus on the transverse measures on the laminations. The map $\varphi$ is proper. Hence we have the following lemma.
  
  \blemm\label{surjectionpropre}\cite[Lem.~10]{Morzy2} The map $\varphi$ defines a continuous, surjective and proper map $\varphi_*$ from the space of Radon measures on
  $[\Gqr]$ to the space of Radon measures 
  on $[\Gmr]$, and then a continuous map $\psi:\mathcal{M}\mathcal{L}_p(\Sigma)\to\mathcal{M}\mathcal{L}_h(\Sigma)$ which is proper and surjective.
  \elemm
  
  A hyperbolic lamination on a compact surface can be endowed with a transverse measure if and only if it has no isolated leaf spiraling on a minimal sublamination.
  Similarly, a measured flat lamination has no isolated leaf terminating in  a minimal sublamination or a cylinder sublamination. Moreover, if a flat lamination $\Lq$ has no isolated leaf which is not periodic, then $\psi(\Lq)$ has no isolated leaf spiraling on a minimal 
  sublamination, then there exists a transverse measure on $\psi(\Lq)$ and hence there exists a transverse measure on $\Lq$. In conclusion, a
  flat lamination on a compact surface can be endowed with a transverse measure if and only if it has no isolated leaf which is not periodic. Notably, 
  every minimal flat lamination of finite graph type $\Lambda$ can be endowed with a transverse measure, and the space of transverse measures on $\Lambda$
  is homeomorphic to the space of transverse measures on the hyperbolic minimal lamination $\psi(\Lambda)$.
  
  \medskip
  
  The group $\RR^{+*}$ acts on $\mathcal{M}\mathcal{L}_p(\Sigma)$ and on $\mathcal{M}\mathcal{L}_h(\Sigma)$ by multiplication of the measures. We denote by  
  $\mathcal{P}\mathcal{M}\mathcal{L}_p(\Sigma)$ and $\mathcal{P}\mathcal{M}\mathcal{L}_h(\Sigma)$
  the quotient spaces for these actions. Since $\psi$ is equivariant by these actions, it defines a continuous map $\overline{\psi}:
  \mathcal{P}\mathcal{M}\mathcal{L}_p(\Sigma)\to
  \mathcal{P}\mathcal{M}\mathcal{L}_h(\Sigma)$ that is surjective and proper. We deduce from this the following lemmas.
  
  \blemm\cite[Lem.~11]{Morzy2} The space $\mathcal{P}\mathcal{M}\mathcal{L}_p(\Sigma)$ is compact.
  \elemm
  
  A {\it measured cylinder lamination} is a measured flat lamination having a unique component which is a cylinder component or a pair of (periodic)
  leaves.
  
  \blemm\cite[Lem.~12]{Morzy2} Since $\Sigma$ is compact, the set of measured cylinder laminations having finitely many leaves is dense in $\mathcal{M}\mathcal{L}_p(\Sigma)$. In particular, 
  $\mathcal{M}\mathcal{L}_p(\Sigma)$ is separable.
  \elemm

  As $\varphi$, the maps $\psi$ and $\overline{\psi}$ are not injective. Assume that $\Sigma$ is compact with genus  $g\in\NN$.
  Considering the lack of injectivity of $\varphi_*$, we see that the preimage of a measured hyperbolic lamination having no closed leaf consists in a unique measured flat lamination
  having no periodic leaf. However, if $(\Lambda_m,\mu_m)$ is a closed leaf $\lambda$ endowed with a transverse measure which, for every arc $\alpha$ such that $G(\alpha)$ 
  contains $\lambda$, is a Dirac measure at $\lambda$ of mass $\delta>0$,
  the preimage of 
  $(\Lambda_m,\mu_m)$ by $\psi$ is the set of measured cylinder laminations whose supports are closed sets $F$ of leaves that are freely homotopic to $\lambda$. 
  If the set of local geodesics of $\srfce$ that are freely homotopic to $\lambda$ contains at least two elements, 
  it foliates a maximal flat cylinder. Then, this set is homeomorphic to $[0,L]$, with $L$ is the height of the flat cylinder, so $F$ is homeomorphic to a closed subset
  of $[0,L]$. Hence, the preimage of $(\Lambda_m,\mu_m)$ by $\psi$ is homeomorphic to the set of Borel measures on $[0,L]$ of mass $\delta$. Since $\varphi_*$ is 
equivariant for the addition of measures, we see that if a measured hyperbolic lamination has some closed leaves $\lambda_1,\dots,\lambda_p$, of respective masses 
 $\delta_1,\dots,\delta_p$  ($p$ is always at most $3g-3$), then the preimage of $(\Lambda_m,\mu_m)$ by $\psi$ is homeomorphic to the Cartesian product of the sets of
 Borel measures
 on $[0,L_i]$, $1\leqslant i\leqslant p$, where $L_i$ is the height of the maximal flat cylinder, union of the local geodesics of $\srfce$ freely homotopic to 
 $\lambda_i$, whose total mass is $\delta_i$.
  
  Since  $\Sigma$ is compact, the projectified space $\mathcal{P}\mathcal{M}\mathcal{L}_h(\Sigma)$ is homeomorphic to the sphere $\SSS^{6g-6}$ 
  (see \cite[Th.~17]{Bonahon97}). 
  If the support of the measured hyperbolic laminations in the equivalence class of $x\in\mathcal{P}\mathcal{M}\mathcal{L}_h(\Sigma)$ has no 
  closed leaf, then the preimage of $\{x\}$
  by $\overline{\psi}$ is a single point. However, if the closed leaves $\lambda_1,\dots\lambda_p$ belong to the support of the measured hyperbolic 
  laminations in the equivalence
  class of $x$, and if $L_1,\dots,L_p$ are the heights of the maximal flat cylinders, unions of the images of the periodic local geodesics of $\srfce$ 
  that are freely homotopic 
  to $\lambda_1,\dots,\lambda_p$ ($L_i=0$ if there is only one periodic local geodesic that is freely homotopic to $\lambda_i$), then the preimage of 
  $\{x\}$ by
  $\overline{\psi}$ is homeomorphic to the Cartesian product with $p$ terms of the sets of Borel measures on  $[0,L_i]$ of total masses at most $1$.
  In particular, since the set of measured hyperbolic laminations whose support contains a closed leaf is projectively dense in the space of measured hyperbolic laminations,
  the subset of $\mathcal{P}\mathcal{M}\mathcal{L}_h(\Sigma)$ of points whose preimages by $\overline{\psi}$ are not single points, is dense in 
  $\mathcal{P}\mathcal{M}\mathcal{L}_h(\Sigma)$.
  
\section{Tree associated to a measured flat lamination.}\label{arbredual}

  Let $\srfce$ be a compact, connected, orientable surface (without boundary), endowed with a half-translation 
structure, such that 
$\chi(\Sigma)<0$, let 
$p:\revet\to\srfce$ be a locally isometric universal cover of covering group $\Gamma_{\widetilde{\Sigma}}$ and let $(\Lambda,\mu)$ be a measured flat lamination on $\srfce$. We denote by 
$(\widetilde{\Lambda},\widetilde{\mu})$ its preimage in $\revet$ and by $\nu_{\widetilde{\mu}}$ the element of $\M_{\grperevet}([\Gqr])$ that it defines
(see the beginning of Section
\ref{liensplatshyperboliques} or \cite[§4]{Morzy2}).
We assume that $\nu_{\widetilde{\mu}}$ has no atom. If it had some, we would blow up them (see the end of \cite[§6]{Morzy2}). 

If $\widetilde{\ell}$ is a leaf of $\widetilde{\Lambda}$, then $\widetilde{\Sigma}-\widetilde{\ell}(\RR)$ has two connected components and the image of a leaf of
$\widetilde{\Lambda}$
is contained in the closure of one of them, since the leaf is not interlaced with $\widetilde{\ell}$. Let $\widetilde{\ell}_0$ and 
$\widetilde{\ell}_1$ be two leaves of $\widetilde{\Lambda}$. We denote by $C_i$ the connected component of $\widetilde{\Sigma}-\widetilde{\ell}_i(\RR)$ that contains 
$\widetilde{\ell}_{i+1}(\RR)$ (with 
$i\in\ZZ/2\ZZ$). We define $C(\widetilde{\ell}_0,\widetilde{\ell}_1)=C_0\cap C_1$. Let $c$ be a geodesic segment joining their images (if
the images are not disjoint, the segment $c$ may be a 
point). A leaf of
$\widetilde{\Lambda}$ {\it intersects $c$ non trivially } if it is contained in $C(\widetilde{\ell}_0,\widetilde{\ell}_1)$ and if
it intersects both complementary components of the image of $c$ in
$C(\widetilde{\ell}_0,\widetilde{\ell}_1)$, and we denote by $B(\widetilde{\ell}_0,\widetilde{\ell}_1)=
B_{\widetilde{\Lambda}}(\widetilde{\ell}_0,\widetilde{\ell}_1)$ the set of leaves of $\widetilde{\Lambda}$ intersecting
$c$ non trivially. 

\blemm The compact set $B(\widetilde{\ell}_0,\widetilde{\ell}_1)$ does not depend on the choice of $c$.
\elemm

If $\{\widetilde{\ell}_1,\widetilde{\ell}_2\}$ is a pair of leaves of $\widetilde{\Lambda}$, we set  
$\widetilde{d}_{\widetilde{\Lambda}}(\widetilde{\ell}_1,\widetilde{\ell}_2)=\frac{1}{2}\nu_{\widetilde{\mu}}(B(\widetilde{\ell}_1,\widetilde{\ell}_2))$. Then 
$\widetilde{d}_{\widetilde{\Lambda}}$ is a pseudo-distance on $\widetilde{\Lambda}$.
We denote by $(T,d_T)$ the quotient metric space $(\widetilde{\Lambda},\widetilde{d}_{\widetilde{\Lambda}})/\sim$, where 
$\widetilde{\ell}\sim\widetilde{\ell}'$ if and only if $\widetilde{d}_{\widetilde{\Lambda}}(\widetilde{\ell},\widetilde{\ell}')=0$ (or equivalently if 
$B(\widetilde{\ell},\widetilde{\ell}')=\{\widetilde{\ell},\widetilde{\ell}'\}$). 

\blemm\label{arbredualaunelamination}\cite[Lem.~16]{Morzy2} The metric space $(T,d_T)$ is an $\RR$-tree.
\elemm

Let $m$ be a hyperbolic metric on $\Sigma$. We recall that $\mathcal{M}\mathcal{L}_h(\Sigma)$ is the space of measured hyperbolic laminations on ${\Sigma}$. 
The set of free homotopy classes of simple closed curves on $\Sigma$ endowed with transverse measures which are Dirac measures of positive masses, embeds into 
$\mathcal{M}\mathcal{L}_h(\Sigma)$, and the intersection number on this set can be extended, in a unique way,  to a continuous map 
$i:\mathcal{M}\mathcal{L}_h(\Sigma)\times\mathcal{M}\mathcal{L}_h(\Sigma)\to\RR^+$ (see \cite[Prop.~3]{Bonahon88}).
  According to Lemma \ref{surjectionpropre}, the map $\varphi_*$ defines a map $\psi:\mathcal{M}\mathcal{L}_p(\Sigma)\to\mathcal{M}\mathcal{L}_h(\Sigma)$.   
   Let $\alpha$ be a non trivial free homotopy class of closed curves, let $(\Lq,\mu_{[q]})$ be a measured flat lamination and let
  $\nu_{\mu_{[q]}}\in\M_{\grperevet}([\Gqr])$ be the measure defined by $(\Lq,\mu_{[q]})$ on $[G_{[\widetilde{q}]}]$.
  We define the intersection number between $(\Lq,\mu_{[q]})$ and $\alpha$ by $$i_{[q]}(\mu_{[q]},\alpha)=
 i(\psi(\Lq,\mu_{[q]}),\alpha).$$ If $k\in\NN$, we have $i_{[q]}(\mu_{[q]},\alpha^k)=k\, i_{[q]}(\mu_{[q]},\alpha)$. Hence, we assume that $\alpha$ is primitive 
 (i.e. if there exists a free homotopy class $\alpha_0$ such that $\alpha=\alpha_0^k$, 
 then $k=\pm1$). We denote by $\alpha_{[q]}$ a periodic local flat geodesic  in the class of $\alpha$, and by $\widetilde{\alpha}_{[q]}$ a lift of
 $\alpha_{[q]}$ in $\widetilde{\Sigma}$. Let $\gamma\in\grperevet-\{e\}$ be one of the two primitive hyperbolic elements of $\grperevet$ whose translation axis is 
 $\widetilde{\alpha}_{[q]}(\RR)$, and let $(\widetilde{\Lambda}_{[q]},\mutilde_{[q]})$ be the preimage of $(\Lambda_{[q]},\mu_{[q]})$ in $\widetilde{\Sigma}$.
  
 \blemm\label{masseintersection}\cite[Lem.~15]{Morzy2} Let $\widetilde{\ell}$ be a leaf of $\Lqr$ which is interlaced with $\widetilde{\alpha}_{[q]}$. 
 The number $i_{[q]}(\mu_{[q]},\alpha)$ is equal to $\frac{1}{2}\nu_{\mu_{[q]}}(B_{\widetilde{\Lambda}}
 (\widetilde{\ell},\gamma\widetilde{\ell})-\gamma\widetilde{\ell})$. If there exists no such leaf, then $i_{[q]}(\mu_{[q]},\alpha)=0$.
 \elemm 
 
  \rem We could define the intersection number between a free homotopy class of closed curves with $(\Lambda_{[q]},\mu_{[q]})$ by the infimum of the masses given
 by the measured
 flat lamination to the closed curves that are piecewise transverse to the lamination, similarly to the intersection number with a measured foliation, 
 but this infimum would not be necessarly attained since the periodic local geodesics are generally not piecewise transverse to the lamination.
 
 Furthermore, in the case of  compact surfaces endowed with a half-translation structure, contrarily to the case of measured hyperbolic lamination (see 
 \cite[Th.~2]{Ota90}), the intersection numbers with the free homotopy classes of closed curves of $\Sigma$ do not separate the measured flat laminations, but only 
 their images in  
 $\mathcal{M}\mathcal{L}_h(\Sigma)$. In particular, the topology defined after Definition \ref{defmesuretransverse} is not equivalent to the one induced by the product topology on 
 $\RR^{{\cal{H}}}$, with ${{\cal{H}}}$ the set of free homotopy classes of closed curves, on the image of $\mathcal{M}\mathcal{L}_p(\Sigma)$ by the map 
 $(\Lambda_{[q]},\mu_{[q]})\mapsto (i(\mu_{[q]},\alpha))_{\alpha\in{{\cal{H}}}}$. 
 
 \medskip
 
The covering group $\Gamma_{\widetilde{\Sigma}}$ acts on $\widetilde{\Sigma}$ by isometries. Hence, it defines an action on the set  
$[\G_{[\widetilde{q}]}]$ of geodesics of $\revet$  that are defined up to changing origin. Since $\widetilde{\Lambda}$ is $\grperevet$-invariant,
this action defines an action on
$\widetilde{\Lambda}$. Since for every $\gamma\in\grperevet$ we have $\gamma_*\nu_{\widetilde{\mu}}=\nu_{\widetilde{\mu}}$ and 
$\gamma B(\widetilde{\ell},\widetilde{\ell}')=B(\gamma\widetilde{\ell},\gamma\widetilde{\ell}')$, for every pair of leaves $\{\widetilde{\ell},\widetilde{\ell}'\}$ of
$\widetilde{\Lambda}$, it defines an isometric action of $\grperevet$ on the  tree $(T,d_T)$ associated to $(\widetilde{\Lambda},\widetilde{\mu})$ defined in Lemma 
\ref{arbredualaunelamination}.

\blemm\cite[Lem.~17]{Morzy2} For every primitive element $\gamma\in\grperevet-\{e\}$, if $\widetilde{\alpha}_{\gamma}(\RR)$ is a tranlation axis of $\gamma$ in $\revet$ and $\alpha_\gamma$
is the projection of $\widetilde{\alpha}_\gamma$ in $\Sigma$, then the translation distance $\ell_T(\gamma)$ of $\gamma$ in $(T,d_T)$ 
is equal to $i_{[q]}(\mu,\alpha_\gamma)$. Moreover, if $\ell_T(\gamma)>0$, the translation axis of $\gamma$ is the image in $T$ of the set of leaves of
$\widetilde{\Lambda}$ which are interlaced with $\widetilde{\alpha}_\gamma$.  
\elemm

If $m$ is a hyperbolic metric on $\Sigma$, and if $(\Lm,\mu_m)$ is a measured hyperbolic lamination on $\srfcem$, then there exists a usual definition of dual tree to 
$(\Lm,\mu_m)$ (see for example \cite[§1]{MorSha91}). The dual tree to $(\Lm,\mu_m)$ can also be defined in a similar way to the tree associated to
a measured flat lamination, and we easily see that the trees we get by the two procedures are isometric by an equivariant isometry (for the isometric actions of the 
covering group, see \cite[§.~8]{Morzy2}).

\medskip

Let $(\Lq,\mu_{[q]})$ be a measured flat lamination on $\srfce$, let $(\Lm,\mu_m)=\psi(\Lq,\mu_{[q]})$, 
and let $(\Lqr,\mutilde_{[\widetilde{q}]})$ and $(\Lmr,\mutilde_m)$ be their preimages in $\widetilde{\Sigma}$. We denote by $\nu_{\mutilde_{[{q}]}}$ and $\nu_{\mutilde_m}$
the elements of $\M_{\grperevet}([\Gqr])$ and $\M_{\grperevet}([\Gmr])$) they define respectivelly. We assume that $\nu_{\mutilde_{[{q}]}}$ and $\nu_{\mutilde_m}$ are 
atomless
(see \cite[§8]{Morzy2} for the general case). The map $\varphi$ defined at section \ref{liensplatshyperboliques}, defines an $\grperevet$-equivariant map 
from $\Lmr$ to $\Lqr$, still denoted by $\varphi$.

By construction, for every pair of leaves $\{\widetilde{\ell}_0,\widetilde{\ell}_1\}$ of $\Lqr$, we have $B_{\Lmr}(\varphi(\widetilde{\ell}_0),
\varphi(\widetilde{\ell}_1))=
\varphi(B_{\Lqr}(\widetilde{\ell}_0,\widetilde{\ell}_1))$. Hence the map $\varphi$ defines a map 
$\varphi_T:(T_{[q]},d_{T_{[q]}})\to(T_{m},d_{T_{m}})$, where 
$(T_{[q]},d_{T_{[q]}})$ and $(T_{m},d_{T_{m}})$ are respectively the tree associated to 
$(\Lqr,\mutilde_{[q]})$ and the dual tree to $(\Lmr,\mutilde_m)$.

\blemm\cite[Lem.~19]{Morzy2} \label{isometrieequivariante}
The map $\varphi_T:(T_{[q]},d_{T_{[q]}})\to(T_{m},d_{T_{m}})$ is a $\grperevet$-equivariant isometry.
\elemm
\bibliographystyle{alphanum}
\bibliography{biblio}{}
Département de mathématique, UMR 8628 CNRS, Université Paris-Sud, Bât. 430, F-91405 Orsay Cedex, FRANCE. 
Bureau : 16.

{\it thomas.morzadec@math.u-psud.fr}
\end{document}